\documentclass[preprint,12pt]{elsarticle}

\usepackage{bm}
\usepackage{amsmath}
\usepackage{amssymb}
\usepackage{mathtools}
\usepackage{siunitx}
\usepackage{fancyref}
\usepackage{pgfplots}
\usepackage{todonotes}
\usepackage{eurosym}
\usepackage{multirow}

\usepackage[acronym]{glossaries}

\usepackage{layout}
\usepackage{printlen}
\usepackage{hyperref}

\usetikzlibrary{math} 

\DeclareSIUnit{\ct}{ct}
\DeclareSIUnit\year{yr}
\DeclareSIUnit{\sieuro}{\text{\euro}}

\journal{Energy}

\newacronym{nlp}{NLP}{Non Linear Problem}
\newacronym{milp}{MILP}{Mixed Integer Linear Program}
\newacronym{minlp}{MINLP}{Mixed Integer Non-Linear Program}
\newacronym{dhn}{DHN}{District Heating Networks}

\begin{document}

\setlength{\parindent}{0cm}

% Specify folder paths

\begin{frontmatter}

\title{Non-linear Topology Optimization of District Heating Networks: A benchmark of Mixed-Integer and Adjoint Approaches}

\author{Yannick Wack \fnref{label1,label2,label3}}
\ead{yannick.wack@kuleuven.be}
\author{Sylvain Serra \fnref{label4}}
\author{Martine Baelmans \fnref{label1,label2}}
\author{Jean-Michel Reneaume \fnref{label4}}
\author{Maarten Blommaert \fnref{label1,label2}}

\fntext[label1]{Department of Mechanical Engineering, KU Leuven, Celestijnenlaan 300 box 2421, 3001 Leuven, Belgium}
\fntext[label4]{Universite de Pau et des Pays de l’Adour, E2S UPPA, LaTEP, Pau, France }
\fntext[label3]{Flemish Institute for Technological Research (VITO), Boeretang 200, 2400 Mol, Belgium}
\fntext[label2]{EnergyVille, Thor Park, Poort Genk 8310, 3600 Genk, Belgium}

%% Abstract
\begin{abstract}
	The widespread use of optimization methods in the design phase of District Heating Networks is currently limited by the availability of scalable optimization approaches that accurately represent the network. In this paper, we compare and benchmark two different approaches to non-linear topology optimization of District Heating Networks in terms of computational cost and optimality gap. The first approach solves a mixed-integer non-linear optimization problem that resolves the binary constraints of pipe routing choices using a combinatorial optimization approach. The second approach solves a relaxed optimization problem using an adjoint optimization approach, and enforces a discrete network topology through penalization. Our benchmark shows that the relaxed penalized problem has a polynomial computational cost scaling, while the combinatorial solution scales exponentially, making it intractable for practical-sized networks. We also evaluate the optimality gap between the two approaches on two different District Heating Network optimization cases. We find that the mixed-integer approach outperforms the adjoint approach on a single-producer case, but the relaxed penalized problem is superior on a multi-producer case. Based on this study, we discuss the importance of initialization strategies for solving the optimal topology and design problem of District Heating Networks as a non-linear optimization problem.
	\end{abstract}

\begin{keyword}
  District Heating Networks\sep  topology optimization \sep mixed-integer non-linear optimization \sep benchmark 
%% keywords here, in the form: keyword \sep keyword
%% MSC codes here, in the form: \MSC code \sep code
%% or \MSC[2008] code \sep code (2000 is the default)
\end{keyword}

\end{frontmatter}

\setcounter{footnote}{0}

% Input variable definition
% variables
% Requiered PACKAGES
%\usepackage{bm}
%\usepackage{mathtools}
%%%%%%%%%%%%%%%%%%%%%%%%%%%%%%%%%%%%%
% Mathematical operators
\newcommand{\tp}[1]{#1^{\intercal}} 			% transpose
\DeclarePairedDelimiter\abs{\lvert}{\rvert} 	% absolute
\newcommand{\card}[1]{\lvert#1\rvert}
\newcommand{\infNorm}[1]{\|#1\|_{\infty}}
\newcommand{\Real}[1]{\mathbb{R}^{#1}}
% Mathematical notation
\newcommand{\ve}[1]{\bm{#1}} 	% vector

%%%%%%%%%%%%%%%%%%%%%%%%%%%%%%%%%%%%%
%% Graph definitions
\newcommand{\gEdge}[3]{#1_{#2#3}}   % generic edge
\newcommand{\gNode}[2]{#1_{#2}}     % generic node

% Directed graph definition
\newcommand{\dirGraph}{G}
\newcommand{\setNodes}{N}
\newcommand{\setEdges}{E}

% Graph components in district heating networks
\newcommand{\pro}{\mathrm{pr}}
\newcommand{\con}{\mathrm{con}}
\newcommand{\rad}{\mathrm{hs}}
\newcommand{\byp}{\mathrm{bp}}
\newcommand{\jun}{\mathrm{jun}}
\newcommand{\pipe}{\mathrm{pipe}}
\newcommand{\operation}{\mathrm{op}}
% subsets
\newcommand{\Npro}{\setNodes_\pro}
\newcommand{\Ncon}{\setNodes_\con}
\newcommand{\NconF}{\setNodes_{\con,\mathrm{f}}}
\newcommand{\NconR}{\setNodes_{\con,\mathrm{r}}}
\newcommand{\Njun}{\setNodes_\jun}
\newcommand{\NproF}{\setNodes_{\pro,\mathrm{f}}}
\newcommand{\NproR}{\setNodes_{\pro,\mathrm{r}}}

\newcommand{\EF}{\setEdges_{\mathrm{f}}}
\newcommand{\Epro}{\setEdges_\pro}
\newcommand{\EproF}{\setEdges_{\pro,\mathrm{f}}}
\newcommand{\EproR}{\setEdges_{\pro,\mathrm{r}}}

\newcommand{\Econ}{\setEdges_\con}
\newcommand{\Erad}{\setEdges_\rad}
\newcommand{\Ebyp}{\setEdges_\byp}
\newcommand{\Epipe}{\setEdges_\pipe}
\newcommand{\EpipeF}{\setEdges_{\pipe,\mathrm{f}}}
\newcommand{\EpipeR}{\setEdges_{\pipe,\mathrm{r}}}
\newcommand{\Eop}{\setEdges_\operation}

% i,j definition
\newcommand{\gi}{i}
\newcommand{\gj}{j}

\newcommand{\giNode}[1]{\gNode{#1}{\gi}}
\newcommand{\gjNode}[1]{\gNode{#1}{\gj}}
\newcommand{\gijEdge}[1]{\gEdge{#1}{\gi}{\gj}}

%%%%%%%%%%%%%%%%%%%%%%%%%%%%%%%%%%%%%%%%%
%% Optimization problem
\newcommand{\cost}{J}
\newcommand{\designVarskal}{\varphi}
\newcommand{\designVar}{\ve{\designVarskal}}
\newcommand{\stateVar}{\ve{x}}
\newcommand{\topVarSkalar}{d}
\newcommand{\topVar}{\ve{\topVarSkalar}}
\newcommand{\pentopVar}{\ve{\penDiameter}}

\newcommand{\equalCon}{\ve{c}}
\newcommand{\inEqualCon}{\ve{h}}

\newcommand{\prodInput}{\gamma} 	% design variable producer inflow
\newcommand{\radValve}{\alpha}    % Heater valve design var
% Reformulation of Optimization problem
\newcommand{\equalNew}{\tilde{\equalCon}}
\newcommand{\inEqualNew}{\tilde{\inEqualCon}}
\newcommand{\designVarNew}{\tilde{\designVar}}

\newcommand{\radFlowNew}{\tilde{\radValve}}
\newcommand{\prodInputNew}{\tilde{\prodInput}}
% Model split up
\newcommand{\equalConModel}{\equalCon_\mathrm{m}}
\newcommand{\inEqualModel}{\equalCon_\mathrm{s}}
\newcommand{\stateVarModel}{\stateVar_\mathrm{m}}
\newcommand{\stateVarIneq}{\stateVar_\mathrm{s}}

% Augmented Lagrangian
\newcommand{\ALagrangian}{\mathcal{L}}
\newcommand{\LagMultis}{\lambda}
\newcommand{\LagPen}{\mu}
\newcommand{\slack}{s}
\newcommand{\equalConState}{g}

%% Cost functions
\newcommand{\costFull}{\mathcal{J}}
\newcommand{\costi}[1]{\cost_{\mathrm{#1}}}

\newcommand{\subCAPEX}{CAP}
\newcommand{\subOPEX}{OP}
\newcommand{\subPipe}{pipe}
\newcommand{\subHeat}{h}
\newcommand{\subPump}{p}
\newcommand{\subRev}{rev}

% NPV
\newcommand{\npv}{NPV}
\newcommand{\npvConst}{f}
\newcommand{\npvN}{A}
\newcommand{\npvt}{t}
\newcommand{\npvDiscount}{e}

\newcommand{\npvConstCAPEX}{\npvConst_{\mathrm{\subCAPEX}}}
\newcommand{\npvConstOPEX}{\npvConst_{\mathrm{\subOPEX}}}

\newcommand{\npvCost}{C}

\newcommand{\actuarialRate}{\npvDiscount_{\mathrm{a}}}
\newcommand{\energyInflation}{\npvDiscount_{\mathrm{i}}}
% Pipe capex
\newcommand{\Jpipepol}{p}
\newcommand{\Jpipesmooth}{k}
\newcommand{\cPipe}{\npvCost_{mathrm{\subPipe}}}

%Heat CAPEX
\newcommand{\cHeatCAPEX}{\npvCost_{\mathrm{hC}}}
\newcommand{\cHeatCAPEXi}[1]{\npvCost_{\mathrm{hC},#1}}
\newcommand{\capacityFactor}{F}
\newcommand{\efficiency}{\eta}
\newcommand{\producerEfficiency}{\efficiency_{\pro}}
%Heat OPEX
\newcommand{\cHeatOPEX}{\npvCost_{\mathrm{hO}}}
\newcommand{\cHeatOPEXi}[1]{\npvCost_{\mathrm{hO},#1}}

% Pump CAPEX and OPEX
\newcommand{\pumpEff}{\efficiency_{\mathrm{pump}}}
\newcommand{\cPumpOPEX}{\npvCost_{\mathrm{pO}}}
\newcommand{\cPumpOPEXi}[1]{\npvCost_{\mathrm{pO},#1}}
\newcommand{\cPumpCAPEX}{\npvCost_{\mathrm{pC}}}
\newcommand{\cPumpCAPEXi}[1]{\npvCost_{\mathrm{pC},#1}}

% Revenue
\newcommand{\cRev}{\npvCost_{\mathrm{r}}}
\newcommand{\cRevi}[1]{\npvCost_{\mathrm{r},#1}}
%%%%%%%%%%%%%%%%%%%%%%%%%%%%%%%%%%%%%%
%% Model equations
\newcommand{\flow}{q}
\newcommand{\pressure}{p}
\newcommand{\temperature}{T}

% Hydraulic properties
\newcommand{\Reynolds}{Re}      % Reynoldsnumber
\newcommand{\density}{\rho}     % water density
\newcommand{\viscosity}{\mu}    % viscosity
\newcommand{\spHeatCap}{c_{\mathrm{p}}}

% Thermal properties
\newcommand{\TOutside}{\temperature_\infty}
\newcommand{\dTinf}{\theta}                 % Difference to outside temperature
\newcommand{\heat}{Q}

% Pipes
\newcommand{\length}{L} % pipe length
\newcommand{\diameter}{d}
\newcommand{\diameterMin}{\diameter_{\mathrm{min}}}
\newcommand{\diameterDiscrete}{D}
\newcommand{\volume}{V}

\newcommand{\rough}{\epsilon} % pipe roughness
\newcommand{\ratioInsul}{r} % insualtion ratio
\newcommand{\condInsul}{\lambda_{\mathrm{i}}} % conductivity thermal insulation
\newcommand{\condGround}{\lambda_{\mathrm{g}}} % conductivity ground

\newcommand{\subScriptOuterD}{o}  % subscript outer pipe diameter
\newcommand{\depthPipe}{h} % depth pipe

\newcommand{\hydrR}{R}  % hydraulic resistance
\newcommand{\frictionfactor}{f} 	% darcy friction factor
\newcommand{\thermR}{U} % thermal resistance

% Pipe junctions
\newcommand{\inflowID}{a}    % Define local index a
\newcommand{\outflowID}{b}

% Heating system
\newcommand{\valveRhydr}{\zeta}   % hydraulic valve resistance

\newcommand{\bypValve}{\alpha}     % bypass valve design variable

\newcommand{\lmtd}{LMTD}
\newcommand{\heaterCoef}{\xi}	% heater coeficiant
\newcommand{\heaterExp}{n} % heater exponent
\newcommand{\THouse}{\dTinf_{\textrm{house}}}

\newcommand{\Qdemand}{Q_{\mathrm{d}}}
\newcommand{\Qdemandi}[1]{Q_{\mathrm{d},#1}}
\newcommand{\DemandSatisfaction}{S}
% Producer
%\newcommand{\prodInput}{q_{\pro}} 	% design variable producer inflow

\newcommand{\prodInputi}[1]{\prodInput_{#1}} 	% design variable producer inflow

\newcommand{\prodTemp}{\Theta}
\newcommand{\prodTempv}{\ve{\Theta}}
\newcommand{\prodTempi}[1]{\Theta_{#1}} % temperature producer inflow

% Variables for benchmark with Merz
\newcommand{\existance}{\phi}
\newcommand{\exPipe}{\phi}
\newcommand{\exProducer}{\omega}

\newcommand{\heatInstalled}{H}
\newcommand{\flowVelocity}{v}
\newcommand{\massFlow}{m}
\newcommand{\bigM}{\mathcal{M}}
\newcommand{\conSimpleHX}{c}

\newcommand{\segment}{s}

% Variables for benchmark
\newcommand{\gams}{fMINLP}
\newcommand{\pathopt}{pNLP}

\newcommand{\walltime}{w}
\newcommand{\walltimeGAMS}{\walltime_{\mathrm{\gams{}}}}
\newcommand{\walltimePATHOPT}{\walltime_{\mathrm{\pathopt{}}}}

\newcommand{\numberPipes}{n}
%%%%%%%%%%%%%%%%%%%%%%%%%%%%%%%%%%%%%%
%% Penalization strategy
\newcommand{\pen}{\xi} % penalization parameter
\newcommand{\penDirection}{a}
\newcommand{\penDiameter}{\bar{\diameter}}
\newcommand{\setDiameters}{\mathcal{S}}
%% Smoothing
\newcommand{\relu}{f}

%Equations
% SI unit definitions

%% SET DEFINITIONS
% Graph subset definition
\newcommand{\defSetPipes}{\forall \gi\gj \in \Epipe}
\newcommand{\defSetRad}{\forall \gi\gj \in \Erad}
\newcommand{\defSetByp}{\forall \gi\gj \in \Ebyp}
\newcommand{\defSetCon}{\forall \gi\gj \in \Erad \cup \Ebyp}
\newcommand{\defSetProd}{\forall \gi\gj \in \Epro}

% Inflow/Outflow definitions
\newcommand{\defInflow}{\inflowID=(\gi,n) \in \setEdges}
\newcommand{\defOutflow}{\outflowID=(n,\gj) \in \setEdges}

%%% Optimization Problem
\newcommand{\topVarBoxConstraints}[1]{\diameterDiscrete_{0}\leq#1\leq\diameterDiscrete_{N}}
\newcommand{\opVarBoxConstraints}[1]{ 0\leq #1 \leq 1}
%%% VARIABLE DEFINITION
% Physical variables
\newcommand{\defFlow}{\ve{\flow} }
\newcommand{\defPressure}{\ve{\pressure} }
\newcommand{\defTemp}{\ve{\dTinf} }

\newcommand{\defStateVar}{\ve{\stateVar} = \left[\ve{\flow},\ve{\pressure},\ve{\dTinf}\right]} 

% Design variables
\newcommand{\defDesignVar}{\designVar=\left[\topVar,\ve{\prodInput}\right]} 
\newcommand{\defTopVar}{\topVar\in\{\diameterDiscrete_{0},\dots,\diameterDiscrete_{N}\}^{\card{\Epipe}}}

% Constraints
\newcommand{\defModelConstraints}{\equalCon\left(\topVar, \designVar,\stateVar\right)}

%% main text
% Introduction
\section{Introduction}

\gls{dhn} are considered one of the core technologies to enable carbon-neutral space heating \cite{OECD/IEA2019}. It has the ability to connect a multitude of different renewable heat sources and provide heat to districts and entire cities. In \gls{dhn}s, the typically high upfront investment cost of groundworks and piping is a core decision variable for the feasibility of a development project. Therefore in the planning phase it is crucial to design the pipe routing (network topology), pipe sizing and heat production capacities in an optimal way. This topology optimization problem solves the question of where to place heating network pipes and at what diameter, while accounting for the non-linear nature of the heat transport physics. Together with the binary choice  of pipe placement it is therefore inherently a \gls{minlp}. 

Solving non-linear combinatorial problems directly can be challenging. As a result, the topology optimization problem of heating networks is often linearized, resulting in a linear optimization problem that can be solved efficiently using Mixed Integer Linear Programming (MILP) solvers.  S\"oderman \cite{Soderman2007}, for example, optimized the structure and configuration of a district cooling network using a MILP approach. Dorfner and Hamacher \cite{dorfner2014large} later used a linear approach to optimize the topology and pipe sizes of a DHN. Haikarainen et al. \cite{Haikarainen2014} optimized the topology and operations, while accounting for different production technologies and heat storage. Mazairac et al. \cite{Mazairac2015a} optimized the topology of a multi-carrier network that incorporates gas and electricity  supply. Morvay et al. \cite{Morvaj2016} optimally designed the network while also optimizing the energy mix supplied. Still using a MILP approach, Bordin et al. \cite{bordin2016optimization} optimized the network topology while studying the set of consumers to optimally connect to a DHN. Most recently, Resimont et al. \cite{Resimont2021} used a MILP approach to optimize a city-scale heating network. The transformation of modern \gls{dhn}s towards multi-source, low-temperature networks \cite{Lund2021} invalidates the assumptions of most linear heating network models. Accurate modeling of heat losses and the different feasible temperature levels of supply and demand sites in \gls{dhn}s requires a non-linear representation of the network physics. Solving the resulting \gls{minlp} problem is challenging. In the literature, one approach to solving this \gls{minlp} problem is to use heuristic search algorithms. For example, Li and Svendsen used a genetic algorithm \cite{Li2013}, and Allen et al. deployed both a minimal spanning tree heuristic and a particle swarm algorithm \cite{Allen2022} to the topology optimization of a small DHN.

To fully leverage the potential of mathematical optimization notably two approaches have been successfully applied to the non-linear optimal topology problem of heating networks. Mertz et al. \cite{mertz2016minlp}, solves the full \gls{minlp} resolving the discrete nature of pipe routing choices, using a combinatorial optimization approach. On the other hand, Blommaert et al. \cite{Blommaert2020a} and Wack et al. \cite{wack2022}, solve a relaxed \gls{nlp} ensuring a discrete network topology through penalization. This method is inspired by partial differential equation (PDE) constrained topology optimization, where penalization methods are combined with adjoint optimization to solve problems governed by PDEs. These methods promise to scale well with the problem size, allowing to optimize large-scale problems.

The current bottleneck to a widespread use of optimization methods for the design and topology optimization of \gls{dhn}s are an accurate representation of the network while maintaining scalability of the approach. The goal of this paper is therefore to compare the performance of the combinatorial approach by Mertz et al. \cite{mertz2016minlp} and the relaxed adjoint approach by Wack et al. \cite{wack2022} to non-linear topology optimization of \gls{dhn}s. First, the computational cost scaling of both optimization approaches with the network size is benchmarked to assess whether they are suited to be used for large-scale \gls{dhn} development projects. Both approaches solve non-linear and non-convex optimization problems, so the found optimal network designs can be different local optima. Therefore, in a second step, the found optimal designs are compared. Here, the influence of resolving the discrete nature of the pipe routing problem, as well as the influence of initialization strategies on the found local optima is discussed. This comparison is performed for a simple single heat producer case and, to closer mimic modern heating network problems, for a multi-producer case with different heat production temperatures.

% Problem definition
\section{The topology optimization problem of \gls{dhn}s}\label{sec:OptProblem}
To compare both approaches, a topology optimization problem is defined that aligns the problem definitions of Mertz et al.\cite{mertz2016minlp} and Wack et al \cite{wack2022}. This ensures comparability of the results with both approaches. For completeness, the optimization problem definition is repeated here in short, for detailed discussion the reader is referred to the afore-mentioned publications.

In order to represent the design and topology optimization problem of a \gls{dhn}, a set of design variables $\defDesignVar$ is defined, containing the pipe diameters $\topVar$ as well as the normalized producer inflows $\ve{\prodInput} $. To represent the physical state of a given network, a vector of physical variables $\defStateVar$ is defined. It contains the flow rates $\defFlow$, nodal pressures $\defPressure$ and nodal and pipe exit temperatures $\defTemp$. The temperatures $\ve{\dTinf} = \ve{\temperature} - \TOutside$ are defined as the difference between the absolute water temperature $\ve{\temperature}$ and the outside air temperature $\TOutside$. Now the topology optimization problem for a \gls{dhn} can be posed as a generic optimization problem of the form:

\begin{equation} \label{eq:opt}
	\begin{aligned} \min_{\designVar,\stateVar}& \qquad
		\costFull \left(\designVar,\stateVar\right) \\
		s.t.& \qquad \inEqualCon(\designVar,\stateVar) \leq 0.  \\
	\end{aligned}
\end{equation}

Here $\costFull$ represents the cost function and is defined as the total cost of the project over an investment horizon of $\npvN = 30 ~\mathrm{years}$. $\inEqualCon(\designVar,\stateVar)$ defines the set of model and technological constraints. A full definition of the cost function and model constraints for this benchmark can be found in \ref{app:model}. This optimization problem is based on Mertz et al. \cite{mertz2016minlp} and Wack et al. \cite{wack2022} and a detailed description of the underlying costs and models can also be found there. In the problem definition, the set definition in table \ref{tab:sets} is used to reference different parts of the network.

\begin{table}
	\centering
	\caption{Notation used to describe different subsets of the \gls{dhn}.}
	\label{tab:sets}
	\begin{tabular}{ll}
		\hline
		Set	& Symbol  \\
		\hline\hline
		All edges &  $\setEdges$ \\  
		Producer input & $\Epro$ \\  
		Pipes& $\Epipe$ \\
		Feed edges &  $\EF$  \\
		Heating system & $\Erad$  \\
		\hline
	\end{tabular}
\end{table}

% Methodology
\section{Methodology}
To compare the two approaches, the topology optimization problem (described in equation \ref{eq:optProblem}) is solved using two different methods: solving the full \gls{minlp} problem in GAMS using the implementation by Mertz et al. \cite{mertz2016minlp} (referred to as \gams{} for simplicity), and using the penalized adjoint optimization approach by Wack et al. \cite{wack2022} (referred to as \pathopt{}). Both implementations were simplified to facilitate comparison. The \pathopt{} implementation no longer accounts for multiple discrete pipe diameters, while the \gams{} implementation no longer optimizes the production side or cascading supply between high and low temperature consumers. The following sections briefly discuss the two methodologies and their different treatment of the binary pipe existence variable.

\subsection*{Combinatorial \gls{minlp} approach: \gams{}}
For the comparison, the topology optimization problem will be solved using a combinatorial \gls{minlp} approach. This approach was previously published by Mertz et al. \cite{mertz2016minlp}, and a comprehensive discussion of the approach can be found there. It is implemented in GAMS and uses the \gls{milp} solver CPLEX and the \gls{nlp} solver CONOPT. To pose the optimization problem as a combinatorial problem, the topological choice of pipe placement is represented by a binary variable $\ve{\exPipe}$. Physical variables defined on these pipes are then coupled to this existence variable using the bigM method ($\bigM$) \cite{mertz2016minlp}, as e.g. in the flow velocity definition: 

\begin{align}\label{eq:optProblem}
	&\gijEdge{\massFlow} -\bigM\left( 1-\gijEdge{\exPipe} \right)\leq\gijEdge{\flowVelocity} \density \pi \frac{1}{4}\gijEdge{\diameter}^2\\ &\gijEdge{\flowVelocity} \density \pi \frac{1}{4}\gijEdge{\diameter}^2\leq \gijEdge{\massFlow} +\bigM\left( 1-\gijEdge{\exPipe} \right) , \quad \defSetPipes.
\end{align}
Here $\bigM \gg \gijEdge{\massFlow}$ is chosen in such a way that if a pipe exists $\gijEdge{\exPipe} = 1$, the bounds on the flow velocity constraint are tight. If $\gijEdge{\exPipe} = 0$, the bounds of this constraints are relaxed. To force the installation of both feed and return pipes the following constraint is used:
\begin{equation}
	\gijEdge{\exPipe}= \gEdge{\exPipe}{\gj}{\gi},\quad\defSetPipes.
\end{equation}

The \gams{} implementation features additional constraints to facilitate convergence. They are defined in a way to not be active in the final optimal design in order to maintain the same final problem formulation as the \pathopt{}. First, the installed capacity $\heatInstalled$ is bound by the sum of consumer demands: 

\begin{equation}
	\sum_{\gi\gj\in\EproF}\gijEdge{\heatInstalled} =\sum_{\gi\gj\in\EproF}\gijEdge{\left(\flow\dTinf  \right)} \density \spHeatCap \leq 1.5 \sum_{kl\in \Erad}\gEdge{Q}{k}{l}.
\end{equation}

To facilitate convergence on the momentum equations, the velocity and pressure drop over a pipe are bound by:

\begin{align}
	\gijEdge{\flowVelocity} &\leq 3.5, \\
	\gijEdge{\flowVelocity} &\leq 18.438 \gijEdge{\diameter} + 0.2186, \\
	\gijEdge{\flowVelocity} &\geq \flowVelocity_{\mathrm{min}}\gijEdge{\exPipe},\\
	\giNode{\pressure} - \gjNode{\pressure}&\leq\SI{200}{\pascal} \quad \defSetPipes.
\end{align}

This constraints are not active in the final optimized design, as they are dominated by the pipe momentum equation (\ref{eq:pipeMomentum}). A fixed pressure drop is set for the heat exchanger

\begin{equation}
	\giNode{\pressure} - \gjNode{\pressure}\geq\SI{2}{\kilo\pascal}, \quad \defSetRad,
\end{equation}

and the exit temperature of the radiator is required to be smaller than the inflow temperature:
\begin{equation}
	\giNode{\dTinf} \geq \gjNode{\dTinf} , \quad \defSetRad.	
\end{equation}

Finally, this optimization strategy uses a sequence of initializations to ensure a stable convergence of the final \gls{minlp} solve, based on the strategy of Marty et al. \cite{marty2018simultaneous}. For completeness, this initialization sequence is briefly described in \ref{sec:Apendix_init}.

\subsection*{Adjoint-based penalization approach: \pathopt{}}

The other considered approach is the adjoint-based penalization approach by Wack et al. \cite{wack2022}. Here the combinatorial problem is relaxed, allowing for a continuous pipe placement choice. A near-discrete topology is then enforced by penalization, e.g. by replacing $\Jpipepol_0$ in the investment cost function (equation \ref{eq:Jpipe}) with
\begin{equation}
	\bar{\Jpipepol}_0(\gijEdge{\diameter})= \Jpipepol_0 \left(\frac{1}{(1+\exp(-\Jpipesmooth\left(\gijEdge{\diameter}-\diameterMin\right))})-1\right)\,,
\end{equation}
or by explicitly penalizing the $\diameter$ in the investment cost function and within the momentum and energy equations as described in Wack et al. \cite{wack2022}. 

Because of this relaxation, the \gls{minlp} reduces to solving a series of \gls{nlp}s. Here, the optimization is initialized from a uniform distribution of pipe diameters.

% Results
\section{Benchmarking computational cost - a single producer case}\label{chap:simple}

Modern \gls{dhn}s grow ever bigger and more complex, including multiple heat production sites at different temperatures, which optimal design tools need to be capable to handle. These tools therefore need to scale well with the heating network size in order to be applicable to relevant cases. In this first benchmark, the computational cost scaling of a combinatorial \gls{minlp} approach (\gams{}) is compared to solving a relaxed penalized \gls{nlp} (\pathopt{}). This is done on a benchmark case containing a single producer. 

\subsection*{Setup}

\begin{figure}[h]
	\includegraphics[width=1\columnwidth]{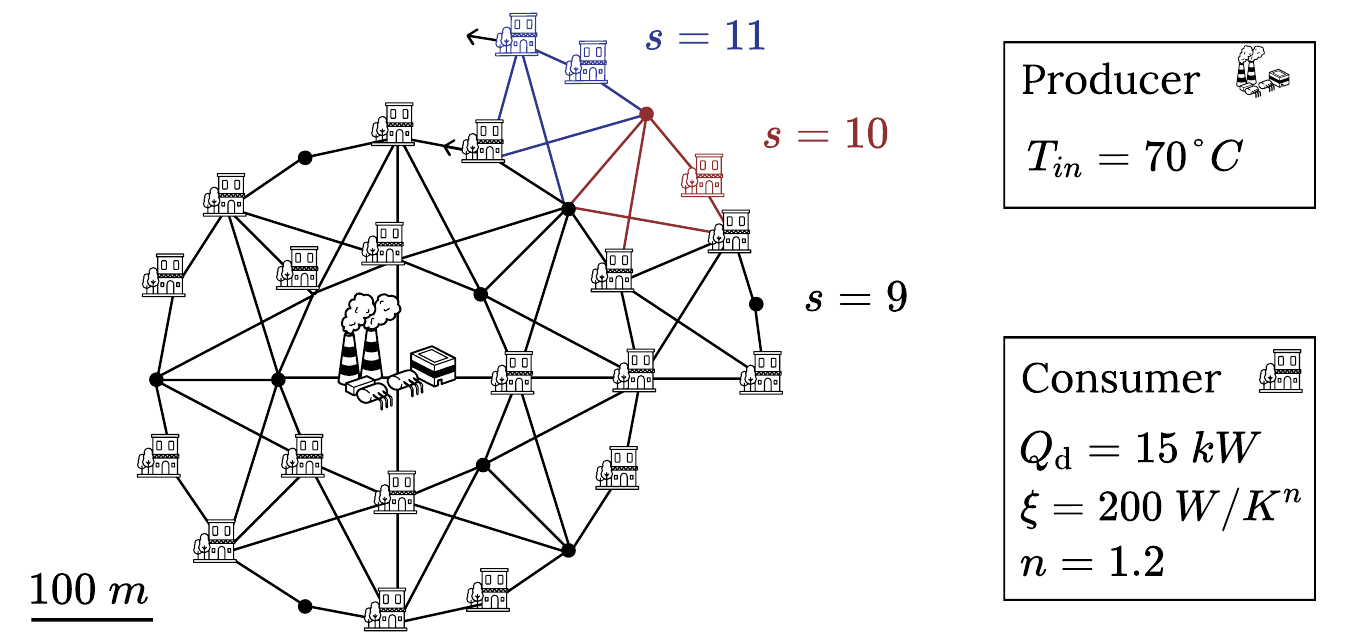}	
	\caption{Setup of the first benchmark case. Around a central producer, heat consumers (houses) are arranged in a circular pattern. They are connected by possible pipe routes (black lines). Pipe junctions are visualized as black circles. To increase the size of this case, additional segments $\segment$ are successively added to the outside of the network.}
	\label{fig:simpleSetup}
\end{figure}

To be able to compare the computational cost scaling of both optimization approaches with increasing network size an easily scalable \gls{dhn} optimization case is set up. This benchmark case is visualized in figure \ref{fig:simpleSetup}. 

Here, a heating network is to be designed around a single heat producer in the center of the network. This producer provides heat at $\SI{70}{\degreeCelsius}$. Around this central producer, houses and possible pipe connections are arranged in a circular manner. All houses have a heat demand of $Q = \SI{15}{\kilo\watt}$, and their heating system characterized by $\heaterCoef=\SI{200}{\watt\per\kelvin^\heaterExp}$ and $\heaterExp=1.2$. To investigate the cost scaling, the size of this circular network is successively increased by adding additional segments $\segment$. With each segment, additional heat consumers and 5 additional potential pipe connections are added, with the number of potential pipes following a linear scaling: $\numberPipes\left( s \right) = 5s+13$. In figure \ref{fig:simpleSetup} the addition of segments $s = 10$ and $s=11$ are visualized in red and blue.

Now the benchmark of both approaches is conducted by creating a series of heating networks of increasing size. Segments $\segment$ are added in steps of 10, starting from $s=0$ up to $s=190$ therefore creating a sequence of heating networks with $(s_k)_{k=0}^{19}, ~s_{k}=10k$ segments. This sequence of optimization problems is solved by both the \gams{} implementation as well as the \pathopt{}. This optimizations were performed on the same computer (Using a single Intel Xeon 3.20 GHz processor core). 

\subsection*{Comparison of the computational cost}
\begin{figure*}[h]
	\centering
	\includegraphics[width=0.736\linewidth]{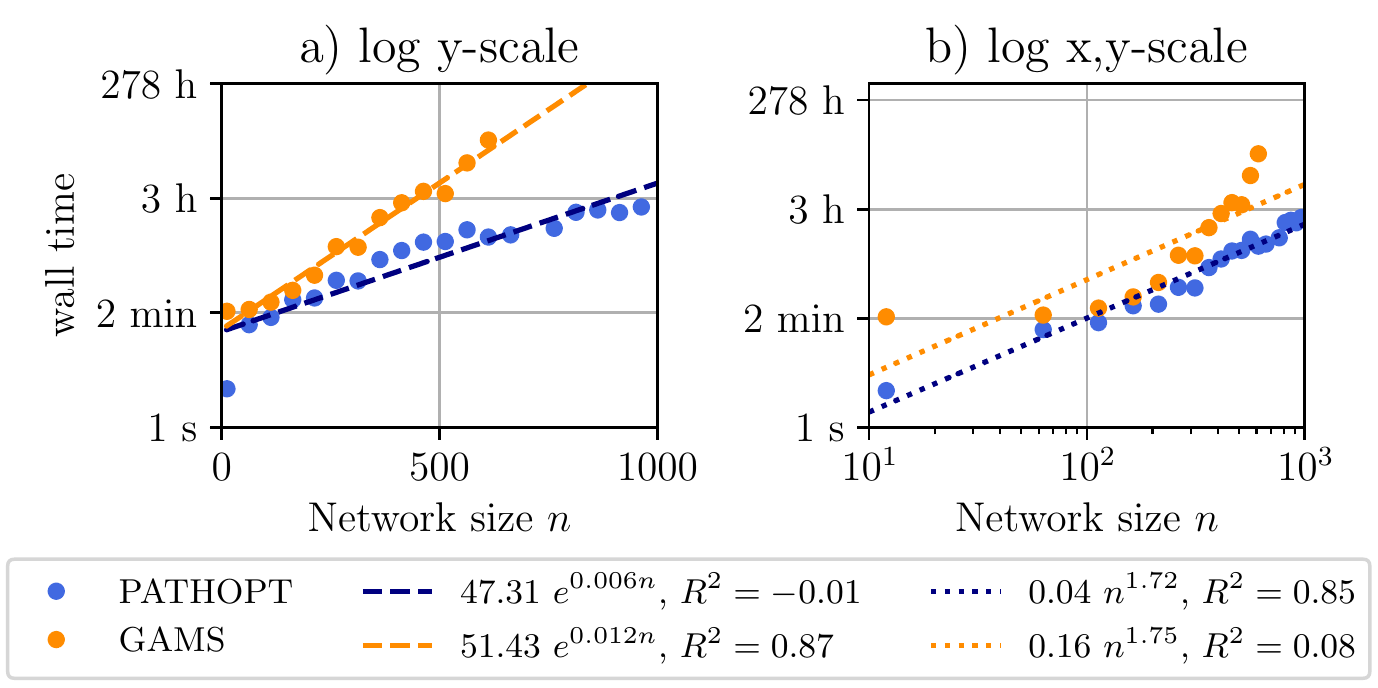}	
	\caption{Wall time scaling of the two approaches with increasing network size $\numberPipes$. The scaling of the approach in the \gams{} is described well by an exponential fit a), while the scaling of the \pathopt{} can be described well by a power fit b). The \pathopt{} did not converge for case $\numberPipes=713$ while the \gams{} did not converge for every case $\numberPipes\geq663$. All non-converged optimizations are excluded from the figure and subsequent analysis}
	\label{fig:simpleScaling}
\end{figure*}

The optimization was repeated 3 times for each network and the mean runtime until convergence (wall time) was recorded. The wall time for each network size is visualized in figure \ref{fig:simpleScaling} on a a) semi-log graph and a b) log-log graph.

The figures show that the wall time for the \gams{} implementation scales faster with the network size then for the \pathopt{}. In order to get an understanding of how long it would take both approaches to optimize networks of considerable size, it is useful to understand what function governs their cost scaling. Considering the algorithmic differences of solving a \gls{minlp} in the case of the \gams{} vs solving a relaxed \gls{nlp} in the case of the \pathopt{}, both an exponential and a polynomial time scaling was tested on both approaches. First, an exponential fit was performed for the scaling of both approaches (see figure \ref{fig:simpleScaling} a)). The wall time scaling of the \gams{} can here be described by the function $\walltimeGAMS\left(  \numberPipes\right) = 51.43e^{0.012\numberPipes} \si{\second}$ with an coefficient of determination $R^2=0.87$, while the wall time scaling of the \pathopt{} can be described by $\walltimePATHOPT\left(  \numberPipes\right) = 47.31e^{0.006\numberPipes}\si{\second}$ with a coefficient of determination of $R^2 = -0.01$. From this fit, it can be concluded that the wall time scaling of the \gams{} approach can be described well by an exponential function, while the wall time scaling of the \pathopt{} approach cannot. Second, a power function was fitted to the wall time scaling of both approaches. Here, the scaling of the \pathopt{} can be described with $\walltimePATHOPT\left(  \numberPipes\right) = 0.04\numberPipes^{1.72}$ with a coefficient of determination of $R^2=0.85$ and the scaling of the \gams{} can be described by $\walltimeGAMS\left(  \numberPipes\right) = 0.16\numberPipes^{1.75}$ with $R^2=0.04$. It can be concluded that the wall time scaling of the \pathopt{} follows a polynomial function reasonably well, while the \gams{} cannot be described with a polynomial function. 

The comparison above highlights the potential exponential solution time of solving full \gls{minlp}s. In the context of \gls{dhn} optimization, this steep cost scaling does not only make design studies slow and expensive, its exponential nature renders optimizations of considerably sized networks intractable and serves as a bottleneck for the optimization of large-scale \gls{dhn}s, containing thousands of potential pipe connections. For example to optimize a network containing $\numberPipes = 2000$ potential pipe connections, an optimization time of $\walltimeGAMS\left(  \numberPipes\right) \approx 43000 ~\mathrm{years}$ would be needed following the observed exponential trend of the \gls{minlp} solution time.
In the context of topology optimization of \gls{dhn}s it is computationally favorable to solve a relaxed \gls{nlp}, as is evident from the polynomial wall time scaling of the \pathopt{} implementation. The polynomial scaling of the computational solution time with this approach makes large-scale \gls{dhn} optimization feasible. This is highlighted for the above mentioned example of a network of $\numberPipes = 2000$, where the necessary optimization time reduces to $\walltimePATHOPT\left(  \numberPipes\right) \approx 5.3 ~\mathrm{hours}$ when using the relaxed \gls{nlp} implementation of the \pathopt{}, assuming the observed polynomial trend.

The potentially exponential computational cost scaling of solving full \gls{minlp}s should be taken into account when choosing optimization strategies for \gls{dhn} topologies. It will make the optimization of large-scale \gls{dhn}s intractable and hinder the application of such optimization tools to real-world \gls{dhn} problems. It is therefore advisable to simplify the \gls{minlp} formulation by either linearizing the non-linear model constraints, posing the problem as a \gls{milp}, or by relaxing the integer constraints, reformulating the problem as a \gls{nlp}. This second relaxation is done within the \pathopt{} implementation, and while keeping the non-linear model in the overall optimization procedure, its significant speed up in the computational cost is shown in this paper.  

\subsection*{Comparison of found optimal network topologies}

\begin{figure}[h]
	\includegraphics[width=1\columnwidth]{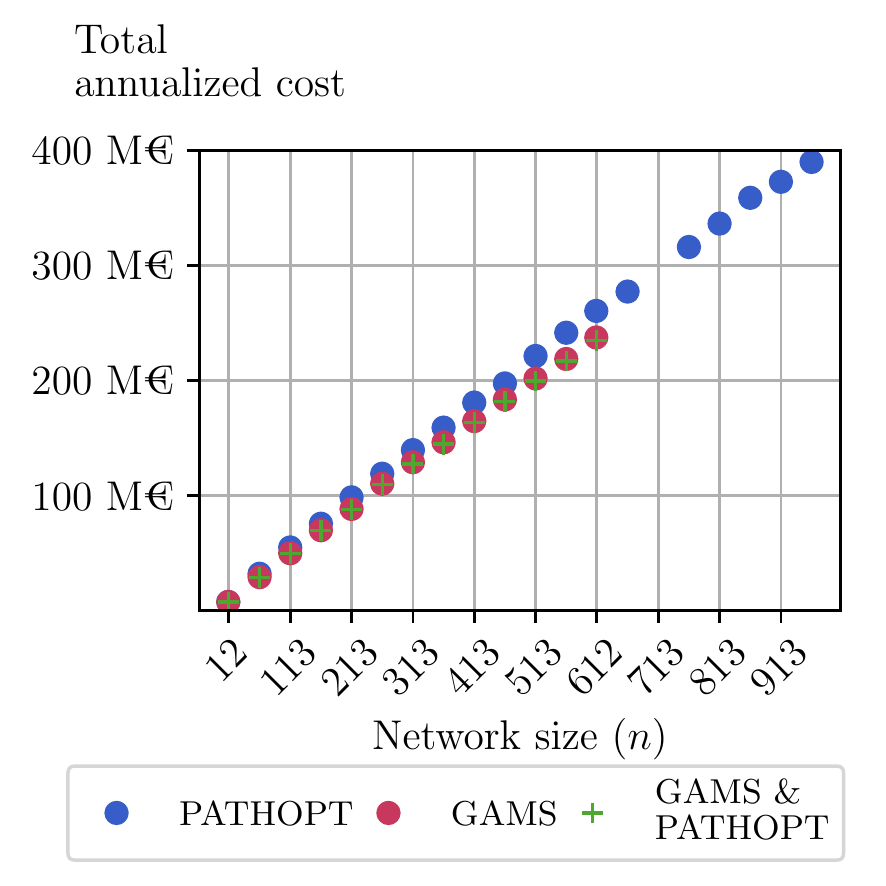}	
	\caption{Comparison of the optimal total annualized cost for each benchmark step found by the \gams{} (orange) and the \pathopt{} implementation (blue). To check if the \gams{} optimum indeed consistently found a cheaper local optimum, the \pathopt{} optimization is run using the \gams{} optima as an initialization (green). The \pathopt{} did not converge for case $\numberPipes=713$ while the \gams{} did not converge for every case $\numberPipes\geq663$. All non-converged optimizations are excluded from the figure and subsequent analysis}
	\label{fig:simpleCost}
\end{figure}

Both approaches solve a non-linear optimization problem. This is necessary because the non-linear heat transport physics can have a big influence on the network design. Non-linearities in the optimization problem, however, can lead to multiple local optima. The proneness of both optimization approaches to these local optima is therefore studied in a next step. Here it can be investigated, if the relaxation of the integer constraints of pipe placement by the \pathopt{} has a relevant influence on the cost of the found optimal network design. For this, the optimality gap of both approaches is compared. The optimal total annualized cost for each benchmark step found with both approaches is visualized in figure \ref{fig:simpleCost}.\footnote{To prevent potential remaining modelling differences from skewing the comparison, both found optimal designs where evaluated in a simulation within the \pathopt{}.}

It can be seen that in this case, the \gams{} implementation (Figure \ref{fig:simpleCost} red) consistently finds a cheaper \gls{dhn} design than the \pathopt{} implementation (Figure \ref{fig:simpleCost} blue). To check if this optimal design found by the \gams{} is indeed a better local optimum, the optimization is repeated in the \pathopt{} and initialized with the optimal design found in the \gams{} (Figure \ref{fig:simpleCost} green). As can be seen in figure \ref{fig:simpleCost}, the \pathopt{} remains in the optimum of the \gams{}, indicating that it is indeed a better local optimum that the \pathopt{} was unable to find with the used initialization strategy. 

\begin{figure*}[h]
	\centering
	\includegraphics[width=0.736\linewidth]{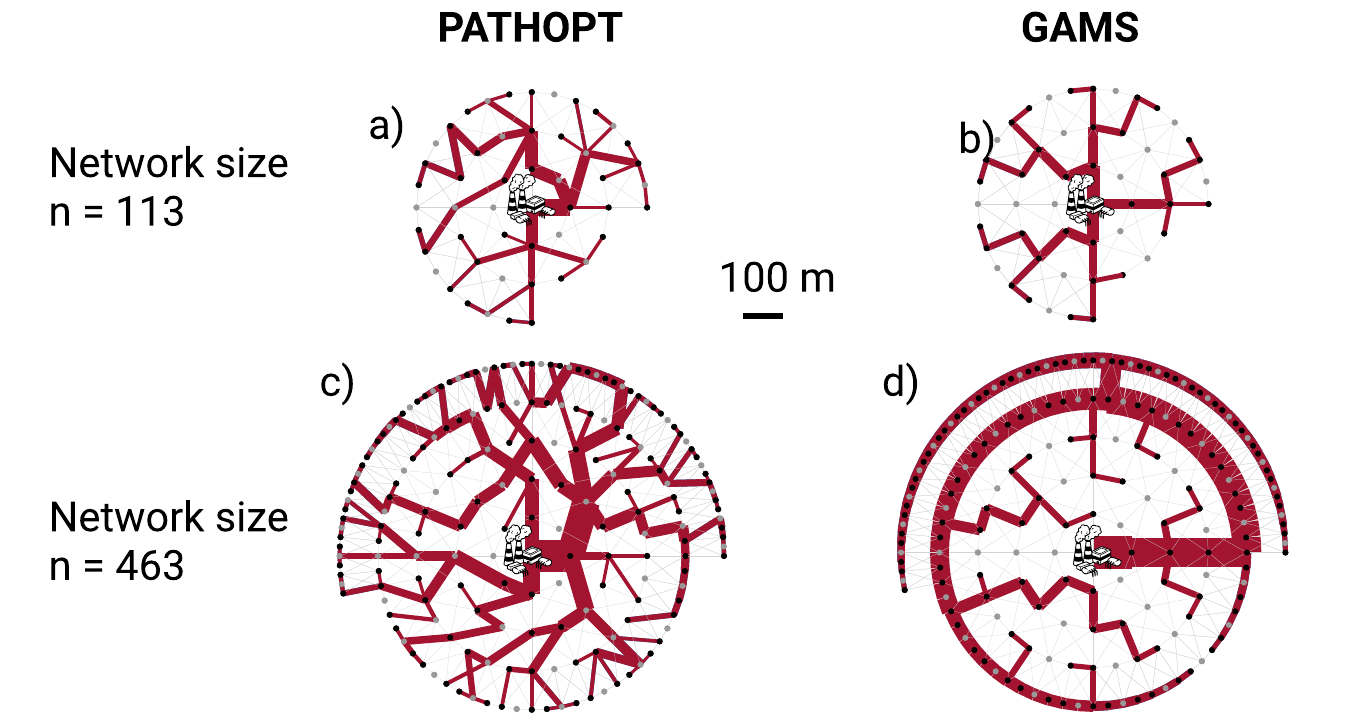}	
	\caption{Comparison of the optimal topology found with the \gams{} and the \pathopt{} for the case of size $\numberPipes = 113$ and the case of size $\numberPipes = 463$. The optimal topologies found by the \gams{} use a shorter pipe length then the topologies found by the \pathopt{}. Red lines represent the pipe network connecting the producer, heat consumers (black dots) and pipe junctions (grey dots). The line thickness indicates the pipe diameter while unused pipe connections are drawn in grey.}
	\label{fig:simpleTopo}
\end{figure*}

To better understand the difference in the optimal network designs causing this cost difference, the optimal network topologies found for network sizes $\numberPipes = 63$ and $\numberPipes = 463$ by both approaches are visualized in figure \ref{fig:simpleTopo}. The figure shows, that the optimal topologies of both approaches differ, with the \pathopt{} using overall a longer pipe length. This difference in the optimal design can be explained with the different initialization methods of the optimization approaches. To facilitate convergence, the \gams{} implementation is initialized with an \gls{milp} solve, that effectively minimizes the overall used pipe length while satisfying the heat demands of the consumers. This initialization promotes topologies with minimal pipe lengths, that turns out to be a good optimal network topology for single producer networks. This optimal topology continues to be an optimum when later solving the full \gls{minlp} as well, as can be seen in figure \ref{fig:simpleTopo} b) and d). The \pathopt{} on the other hand, initializing the full \gls{minlp} from a uniform distribution of pipes, finds a different local optimum using an overall larger pipe length (see figure \ref{fig:simpleTopo}). This knowledge of the influence of initializations on the optimal topology, can be used to reduce the proneness to local optima for both optimization approaches. In the following section \ref{sec:2producer} it will be studied how the initialization of both approaches perform on a multi producer case.

\section{The influence of initialization - a two producer case}\label{sec:2producer}
Modern $4^{\text{th}}$ generation \gls{dhn}s often feature multiple producers with different injection temperatures. To compare how both optimization approaches perform when designing such networks, a benchmark case was designed featuring two heat producers. 

\subsection*{Setup}
\begin{figure}[h]
	\includegraphics[width=1\columnwidth]{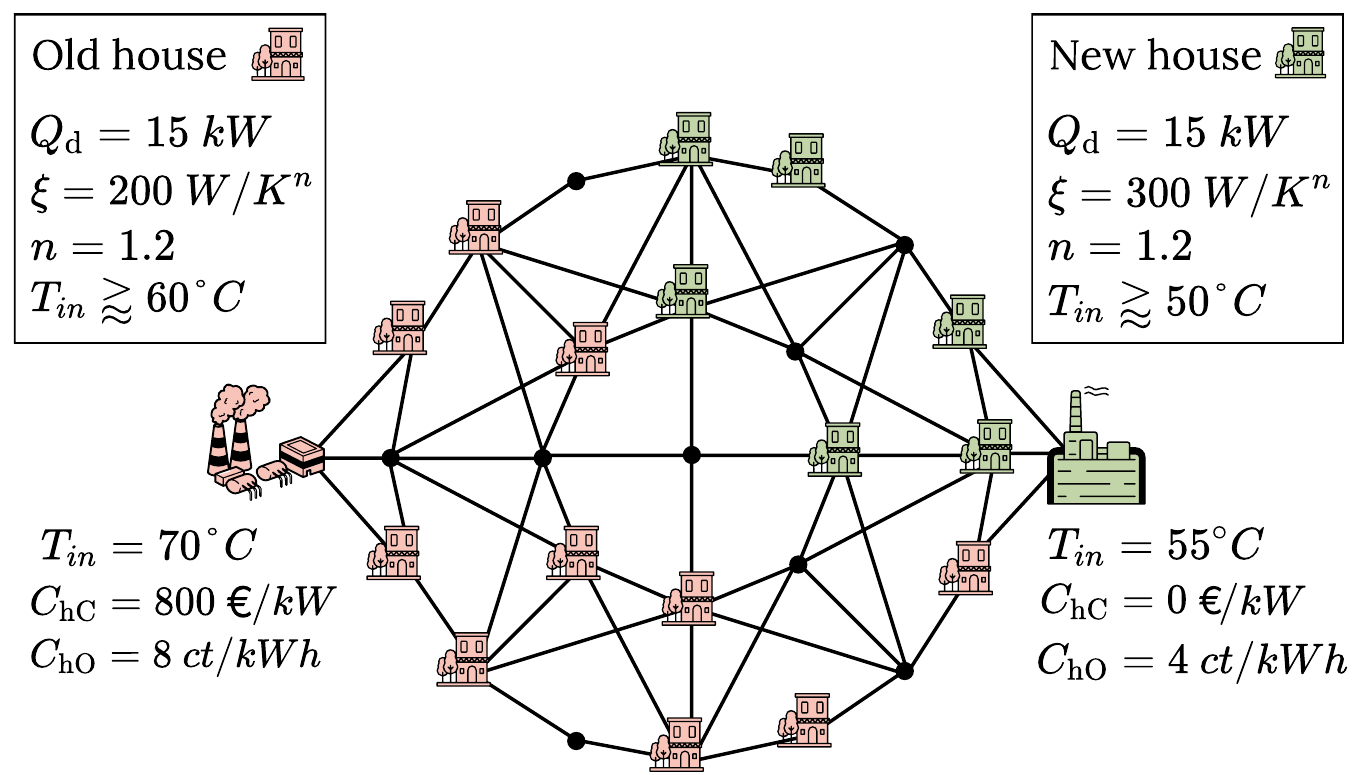}	
	\caption{Setup for the two producer case. This setup features two producers on the left and right of a circular network. The network features old houses that require high temperature heat (red) and newer houses that can satisfy their demand with lower temperatures (green). The heat from the left producer (red) can supply all houses, while the right producer (green) can only supply the houses in the top right quadrant}
	\label{fig:hardSetup}
\end{figure}

For this benchmark, a circular scalable test case with two producers with different injection temperatures is used. The layout of this case is shown in figure \ref{fig:hardSetup}. The circular network superstructure for this case is equivalent to the case in section \ref{chap:simple}, representing a potential district to be connected to a \gls{dhn}. Here, two producers are placed on the left and right of this district. The two producers represent different heat production technologies, with the left producer supplying heat at $\SI{70}{\degreeCelsius}$ at high costs ($\cHeatCAPEX = \SI[per-mode=symbol]{800}{\sieuro\per\kilo\watt}$, $\cHeatOPEX = \SI[per-mode=symbol]{8}{\ct\per\kilo\watt\per\hour}$) while the producer on the right supplies heat at $\SI{55}{\degreeCelsius}$ at a lower cost ($\cHeatCAPEX = \SI[per-mode=symbol]{0}{\sieuro\per\kilo\watt}$, $\cHeatOPEX = \SI[per-mode=symbol]{4}{\ct\per\kilo\watt\per\hour}$). The district contains houses with different heat system characteristics. While the houses on the top right, using a modern heating systems, can work with water at $\gtrapprox\SI{50}{\degreeCelsius}$, the rest of the houses require heat at higher temperatures ($\gtrapprox\SI{60}{\degreeCelsius}$) than the low temperature producer can provide. The size of this network can be increased by adding additional rings of houses to this network, representing different sizes of DHN development projects.

\subsection*{Comparison of the optimized network designs}
\begin{figure*}[h]
	\centering
	\includegraphics[width=0.736\linewidth]{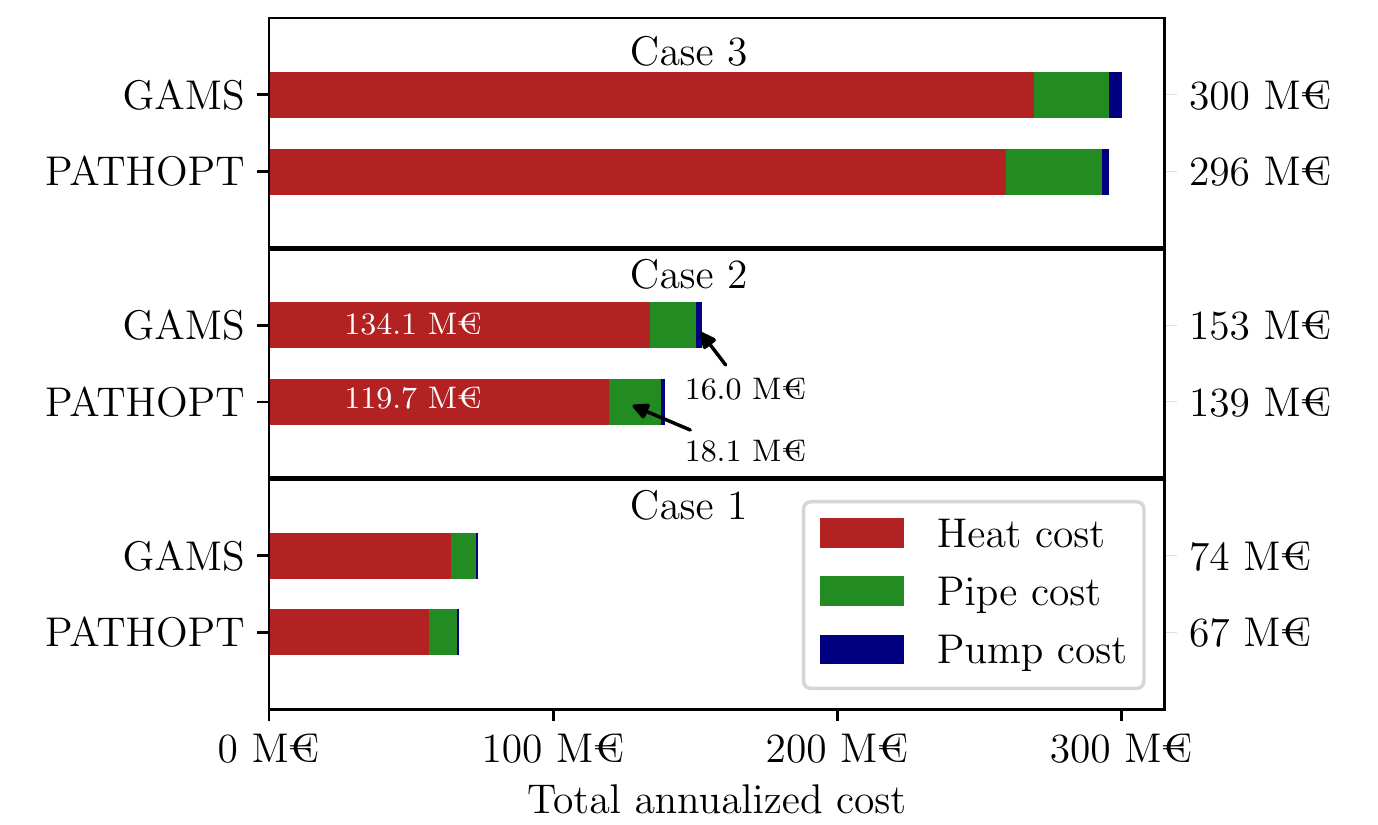}	
	\caption{Comparison of the optimal total annualized cost for case 1-3 found by the \gams{} and the \pathopt{} implementation. While the \gams{} finds optimal design of lower pipe cost (green), the \pathopt{} finds optimal design that save on heating costs (red), ultimately leading to a lower total annualized cost for all three cases.}
	\label{fig:hardCost}
\end{figure*}

Now this case is optimized using both the \gams{} and the \pathopt{} implementation. To compare the two approaches, three different networks of increasing sizes where optimized (Case 1 with $\numberPipes=138$, case 2 with $\numberPipes = 298$ and case 3 with $\numberPipes = 618$). The annualized costs of the resulting optima found by both approaches are visualized in figure \ref{fig:hardCost}.\footnote{Again the optima of both approaches where evaluated with a \pathopt{} simulation to avoid discrepancies.}

The cost comparison shows that in this two producer case, the \pathopt{} found an optimal network design with a lower total annualized cost then the \gams{} in all three cases. While the optimal network design of the \gams{} is again cheaper in pipe investment cost, large savings in heat costs can be made with the proposed optimal design by the \pathopt{}. While in case 2 the pipe investment cost of the optimal network design found by the \gams{} is $\SI{2.1}{\mega \sieuro}$ lower, the optimal design found by the \pathopt{} achieves heat cost savings of $\SI{14.4}{\mega \sieuro}$, ultimately leading to a lower total annualized cost. 

\begin{figure*}[h]
	\centering
	\includegraphics[width=0.736\linewidth]{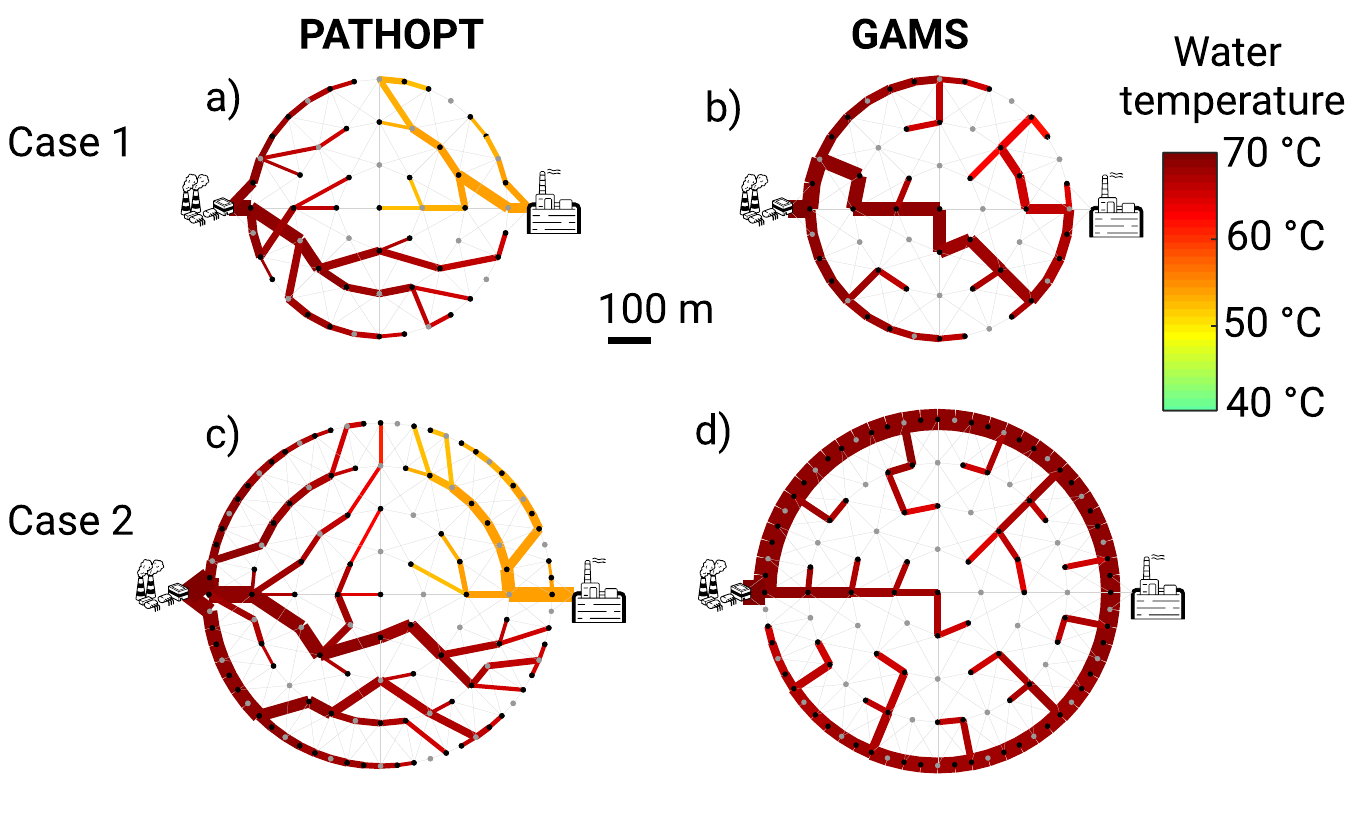}	
	\caption{Colored lines represent the pipe network connecting the producer, heat consumers (black dots) and pipe junctions (grey dots). The line thickness indicates the pipe diameter while unused pipe connections are drawn in grey. The line color corresponds to the water temperature in the pipe. While the \pathopt{} finds an optimal topology connecting modern houses to the low temperature source, the \gams{} implementation favors a single network connecting all houses to the high temperature source.}
	\label{fig:hardTopo}
\end{figure*}

For a better understanding of the difference in the optimal designs, the optimal topologies of both approaches are compared. In figure \ref{fig:hardTopo}, the optimal network topologies of case 1 and 2 optimized by both the \pathopt{} and the \gams{} are visualized. It can be seen that the optimal topologies found by the \pathopt{} connect most new houses to the low temperature producer on the top right, effectively creating two separate networks at different temperatures. This use of cheaper, low temperature heat from the producer on the right leads to the cheaper optimal network designs in comparison to the \gams{} that where observed in figure \ref{fig:hardCost}. The initialization steps of the \gams{}\footnote{Further described in \ref{sec:Apendix_init}} on the other hand, minimizes pipe length while satisfying the heat demands, therefore favoring an optimal topology for both cases that connects all houses to the high temperature producer on the left. The topology resulting from this initialization remains a local optimum also for the full \gls{minlp} optimization. While saving on pipe investment cost, this optimal network topology found by the \gams{}, is more expensive than the optimal topology by the \pathopt{} because it relies on buying more high temperature expensive heat from the producer on the left.   

The results of the benchmarking study indicate that the initialization strategy for the topology optimization problem of \gls{dhn}s has a significant impact on the optimal designs found for both the single- and two-producer cases. This dependence on the initialization should always be taken into account when solving non-convex optimal design problems for DHNs that use non-linear models or cost functions. This sensitivity to the initialization could be mitigated by using globalization strategies, such as running multiple optimizations from different initializations, as demonstrated by Marty et al. \cite{marty2018simultaneous}.

 % Conclusion
 \section{Conclusion}
 In this paper, we compared the performance of two different approaches to the non-linear topology optimization of \gls{dhn}s. While the approach by Mertz et al. \cite{mertz2016minlp} solves a full \gls{minlp}, resolving the discrete nature of pipe routing choices using a combinatorial optimization approach (\gams{}), the approach by Wack et al \cite{wack2022} solves a relaxed \gls{nlp} ensuring a discrete network topology through penalization (\pathopt{}). By performing a benchmark of both approaches, we showed that in the context of topology optimization of \gls{dhn}'s, solving the full \gls{minlp} leads to an exponential scaling of the computational cost with the network size. This steep cost scaling does not only make design studies slow and expensive, its exponential nature renders optimizations of large networks containing thousands of potential pipe connections intractable. On the other hand, solving the topology optimization problem as a relaxed \gls{nlp}, as is done by the \pathopt{}, maintains a polynomial scaling of the computational cost. This polynomial scaling of the optimization time makes large-scale \gls{dhn} optimization feasible. We highlighted the scaling difference of both approaches by showing that for a \gls{dhn} featuring $\numberPipes=2000$, an optimization time of $\walltimeGAMS\left(  \numberPipes\right) \approx 43000 ~\mathrm{years}$ would be needed following the observed exponential trend of the \gls{minlp} solution time in the \gams{} implementation. The necessary optimization time reduces to $\walltimePATHOPT\left(  \numberPipes\right) \approx 5.3 ~\mathrm{hours}$ when using the relaxed \gls{nlp} implementation of the \pathopt{}, assuming the observed polynomial trend.
 
 Posing the optimization problem as a non-linear optimization problem has the benefit of an accurate representation of the \gls{dhn} physics. Non-linearities in the optimization problem can lead to multiple local optima. We showed in this paper how the initialization of different optimization approaches can therefore lead to different optimal designs. For a \gls{dhn} optimization case containing a single producer, the \gams{} implementation consistently found cheaper optimal network designs than the \pathopt{} approach. This difference in optimal design was explained by the different initialization strategies, notably the \gams{} implementation focusing on minimizing pipe length in the initialization. The resulting optimal topologies featuring minimal pipe lengths, prove beneficial for single producer networks. For a two producer case on the other hand, this initialization strategy showed to be a limitation and the \pathopt{} approach consistently found cheaper optimal network topologies than the \gams{} implementation. The comparison of the optimality gap between both approaches shows that resolving the binary nature of the pipe routing problem by solving a full \gls{minlp} does not necessarily lead to better optimal \gls{dhn} topologies, while significantly increasing the computational cost. 
 
 The benchmark also showed that the proneness to local optima and the consequent importance of initialization strategies should be taken into account when deciding to solve the optimal design problem of \gls{dhn}s as a \gls{nlp}. This proneness to local optima can be mitigated with globalization strategies like e.g. running multiple optimizations from different initializations. Although used in literature, in this paper heuristic algorithms and linearized models were omitted. Future work should complement this comparison by benchmarking the scalability of said heuristics (e.g. genetic algorithms) and \gls{milp} approaches to the topology optimization of \gls{dhn}'s. 
 
 % Aknowledgements
 \section*{Data Availability}
 A data-set including the structure, input parameters and optimization results of the heating networks used in both benchmark cases of this paper is available at the following link: \url{https://doi.org/10.48804/DO1BRQ}. The optimization results can be replicated using the methodology and formulations described in this paper.
 
 Additionally a small case generator was written that can be used to recreate all benchmark cases of this paper. The repository can be found at the following link: \url{https://doi.org/10.5281/zenodo.7434451}

 \section*{Acknowledgements}
A special thank you to Valentin Gressot for his support in aligning the software base of both approaches for the benchmark.

Yannick Wack is funded by the Flemish institute for technological research
(VITO). 

\section*{CRediT authorship contribution statement}
\textbf{Yannick Wack:} Conceptualization, Formal analysis, Software, Visualization,
Writing – original draft. \textbf{Sylvain Serra:} Conceptualization, Formal analysis, Software,
Writing - Review \& Editing. \textbf{Martine Baelmans:} Conceptualization, Writing - Review \& Editing.  \textbf{Jean-Michel Reneaume:} Conceptualization, Writing - Review \& Editing.  \textbf{Maarten Blommaert:} Conceptualization, Software, Writing - Review \& Editing, Supervision. 

\section*{Compliance with ethical standards}
\subsubsection*{Conflict of interest}
The authors declare that they have no conflict of interest
\subsubsection*{Funding}
The authors did not receive support from any organization for the submitted
work.
\subsubsection*{Ethical approval}
This article does not contain any studies with human participants or animals
performed by any of the authors.

%% The Appendices part is started with the command \appendix;
%% appendix sections are then done as normal sections
\appendix
\section{Detailed definition of the topology optimization problem for \gls{dhn}s}\label{app:model}
In this section, the underlying cost function $\costFull$ and the set of model and technological constraints $\inEqualCon(\designVar,\stateVar)$ for this benchmark are defined in detail. 

\subsection{Cost function}\label{sec:Cost}
The objective function is defined as the total cost of the project over an investment horizon of $\npvN = 30 ~\mathrm{years}$:
\begin{align} \label{eq:totalCost}
	\costFull\left(\designVar,\stateVar\right) &= \npvConstCAPEX \left(\costi{\subPipe,\subCAPEX}\left(\designVar\right) + \costi{\subHeat,\subCAPEX}\left(\stateVar\right)  \right)\nonumber\\
	&+ \npvConstOPEX \left(\costi{\subHeat,\subOPEX}\left(\stateVar\right) + \costi{\subPump,\subOPEX}\left(\stateVar\right)\right) \,,
\end{align}

with

\begin{align} \label{eq:npvDiscount}
	\npvConstCAPEX =& \left(1+\actuarialRate\right)^\npvN ,\\ \npvConstOPEX =& \frac{1-\left(1+\actuarialRate\right)^\npvN\left(1+\energyInflation\right)^\npvN}{1-\left(1+\actuarialRate\right)\left(1+\energyInflation\right)}  \, ,
\end{align}

assuming a discount rate $\actuarialRate = 0.04$, and an energy inflation rate $\energyInflation =0.04$ \cite{mertz2016minlp}. The invest cost of piping $\costi{\subPipe,\subCAPEX}$ is approximated with a linear interpolation of the catalogue cost of commercially available pipe diameters and the trench cost:  
\begin{equation} \label{eq:Jpipe}
	\costi{\subPipe,\subCAPEX}\left(\designVar\right) = \sum_{\gi\gj \in \Epipe}\left(\Jpipepol_1  \gijEdge{\diameter} + \Jpipepol_0\right)\gijEdge{\length}\,. 
\end{equation}
With the interpolation coefficients $\Jpipepol_1= \SI{1976.3}{\sieuro\per\meter\squared}$, $\Jpipepol_0=\SI{301.4}{\sieuro\per\meter}$. The investment cost for building heat production plants is calculated using %$\Jpipepol_0=\SI{301.4}{\sieuro\per\meter}$ and $\Jpipepol_{\text{tr}}=\SI{400}{\sieuro\per\meter}$
\begin{equation}
	\costi{\subHeat,\subCAPEX}\left(\stateVar\right) = \frac{\density \spHeatCap}{\capacityFactor \producerEfficiency}  \sum_{\gi\gj \in \EproF}\left( \cHeatCAPEX\,\flow \,\Delta\dTinf \right)_{\gi\gj}   \, ,
\end{equation}
with the capacity price of heat production $\cHeatCAPEX = \SI{800}{\sieuro\per{\kilo\watt}}$, the capacity factor $\capacityFactor=0.33$ and assuming an efficiency of $\producerEfficiency = 0.9$. The operational heat cost is calculated using
\begin{equation}
	\costi{\subHeat,\subOPEX}\left(\stateVar\right) = \frac{\density \spHeatCap}{\producerEfficiency}8760\frac{\si{\hour}}{\si{\year}} \sum_{\gi\gj \in \EproF}\gijEdge{\left(\cHeatOPEX \,\flow \, \Delta\dTinf  \right)}\, ,
\end{equation}
with the unit price of heat $\cHeatOPEX =\SI{0.06}{\sieuro\per{\kilo\watt\hour}} $. The operational cost of pumps at the heat production sites is computed with
\begin{equation}
	\costi{\subPump,\subOPEX}\left(\stateVar\right) =  \frac{8760}{\pumpEff}\frac{\si{\hour}}{\si{\year}}\sum_{\gi\gj \in \EproF}\cPumpOPEXi{\gi\gj}\left(\gjNode{\pressure}-\pressure_{a}\right) \gijEdge{\flow} \, ,
\end{equation}
with the electricity price $\cPumpOPEX=\SI{0.11}{\sieuro\per{\kilo\watt\hour}} $ and a pump efficiency of $\pumpEff = 0.7$.

\subsection{\gls{dhn} model constraints}\label{sec:Model}
For all pipe junctions in the network conservation of mass and energy is assumed.
\subsubsection*{Pipe model}
The momentum equations in the pipes are modelled using the Blasius friction factor $\frictionfactor$ and assume a singular pressure drop of $30\%$ \cite{mertz2016minlp}:
\begin{equation} \label{eq:pipeMomentum}
	(\giNode{\pressure} - \gjNode{\pressure}) =  \frac{100}{70}\gijEdge{\frictionfactor} \frac{8\density\gijEdge{\length}}{\gijEdge{\diameter}^5\pi^2}\abs{\gijEdge{\flow}}\gijEdge{\flow} \,,\quad \defSetPipes,
\end{equation}
\begin{equation}
	\mathrm{with} \quad	\gijEdge{\frictionfactor} = 0.3164\left(\Reynolds\right)^{-\frac{1}{4}} \,,\quad \defSetPipes.
\end{equation}
The energy equation over pipes are modelled accounting for the thermal conductivity of the pipe insulation $\condInsul = \SI{0.03}{\watt\per\meter\per\kelvin}$ and the surrounding soil $\condGround = \SI{1.4}{\watt\per\meter\per\kelvin}$,
\begin{equation}\label{eq:pipeEnergy}
	\gijEdge{\dTinf} = \giNode{\dTinf} \exp{\left(\frac{-\gijEdge{\length}}{\density \spHeatCap \abs{\gijEdge{\flow}}{\gijEdge{\thermR}}}\right)}, \quad \defSetPipes,
\end{equation}
\begin{equation}
	\gijEdge{\thermR} = \frac{\ln(4 \depthPipe/(\ratioInsul \gijEdge{\diameter}))}{2 \pi \condGround} +\frac{\ln{\ratioInsul}}{2 \pi \condInsul},
\end{equation}
assuming an insulation ratio $\ratioInsul = 1.4$ and a pipe depth $\depthPipe = \SI{0.4}{m}$ \cite{Blommaert2020a}.
\subsubsection*{Consumer model}
In the consumer arc, conservation of momentum is assumed. Conservation of energy in the heating system leads to
\begin{equation}\label{eq:heaterConservation}
	\density \spHeatCap \gijEdge{\flow}(\giNode{\dTinf} - \gijEdge{\dTinf}) = \gijEdge{\heat}, \quad \defSetRad,
\end{equation} 
with $\gijEdge{\heat}$ the heat transferred to the house through the heating system.
The latter is modelled with the characteristic equation for radiators \cite{Ashrae2014c} using the $\lmtd$ approximation by Chen \cite{Chen1987}:
\begin{equation}\label{eq:characteristicRadiator}
	\gijEdge{\heat} = \gijEdge{\heaterCoef} \left(\lmtd\left(\giNode{\dTinf} - \THouse,\gijEdge{\dTinf} - \THouse\right)\right)^{\gijEdge{\heaterExp}}.
\end{equation}
\begin{align}
	\mathrm{with}\quad&\lmtd\left(\Delta\dTinf_{\mathrm{A}},\Delta\dTinf_{\mathrm{B}}\right)\\ &\approx\left(\Delta\dTinf_{\mathrm{A}}\Delta\dTinf_{\mathrm{B}}\left(\frac{\Delta\dTinf_{\mathrm{A}}+\Delta\dTinf_{\mathrm{B}}}{2}\right)\right)^{\frac{1}{3}}\,.
\end{align} 

Here, $\THouse$ is the inside temperature of the house, and $\heaterCoef$ and $\heaterExp$ are heating system specific coefficients.

\subsubsection*{Producer model}
At the producer, heat is injected with an input flow $\prodInput$  at a fixed temperature $\prodTemp$:
\begin{equation} \label{eq:12}
	\gijEdge{\flow} = \prodInputi{\gi\gj},\quad \gijEdge{\dTinf} = \prodTempi{ij} \quad \forall \gi\gj \in 
	\EproF.
\end{equation}

\subsection*{Additional state constraints}\label{sec:stateConst}
To ensure that the heat demand $\Qdemandi{\gi\gj}$ of every consumer is met, the following constraint is defined:
\begin{align}\label{eq:heatSatisfaction}
	\gijEdge{Q}-\Qdemandi{\gi\gj} \geq 0 \,, \quad \defSetRad.
\end{align}

\section{Initialization strategy of the combinatorial approach \\(\gams{})} \label{sec:Apendix_init}
\begin{figure}[h]
	\centering
	\includegraphics[width=\columnwidth]{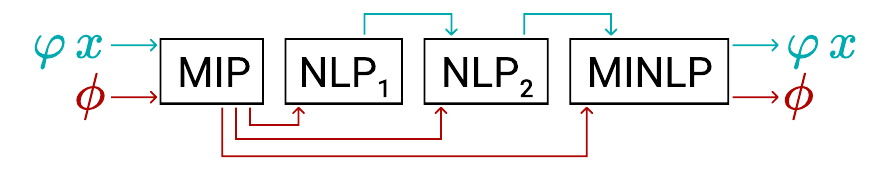}
	\caption{Initialization strategy of the combinatorial approach (\gams{}). The optimal network topology $\existance$ determined by the MIP is used to initialize the following NLPs (red). The design variables $\designVar$ and physical variables $\stateVar$ are optimized in a sequence of NLPS (blue).}
	\label{fig:initializationStrategy}
\end{figure}

The initialization strategy for the \gls{minlp} approach is crucial for the overall solution process. The initial values for the variables and constraints in the MINLP model can greatly affect the solution's speed, accuracy, and feasibility. Figure \ref{fig:initializationStrategy} shows the specific initialization strategy used in the \gams{} implementation.

The initialization strategy for the combinatorial \gams{} implementation consists of the following steps. In the first step (MIP), a \gls{milp} model is solved to determine an initial network topology. For simplicity, the non-linearities present in the system are not considered. The optimal network topology of the MIP step is used in $\mathrm{NLP_1}$,$\mathrm{NLP_2}$ and the MINLP as the initial network topology (Visualized by the red arrows in figure \ref{fig:initializationStrategy}). The first NLP model ($\mathrm{NLP_1}$) is designed to consider some of the non-linearities that were ignored in the initial MIP model. In this step, the characteristic equation of a radiator is used in a simplified form in order to reduce the computational complexity of the optimization problem. The second NLP model ($\mathrm{NLP_2}$) takes into account the full set of non-linearities present in the system. Finally, the full \gls{minlp} is solved including all constraints and binary variables.

%% References
%%
%% Following citation commands can be used in the body text:
%% Usage of \cite is as follows:
%%   \cite{key}          ==>>  [#]
%%   \cite[chap. 2]{key} ==>>  [#, chap. 2]
%%   \citet{key}         ==>>  Author [#]

%% References with bibTeX database:

%%-----------------------------
%%      your bibliography
%%-----------------------------
\bibliographystyle{ieeetr}
\bibliography{library}% common bib file

\end{document}